\documentclass[11pt]{amsart}
\pdfoutput=1
\usepackage{comment}
\usepackage{amsmath, amsfonts, amssymb}
\usepackage{cellspace}
\usepackage{svg}
\usepackage{csquotes}
\usepackage[colorlinks=true, citecolor=cyan, urlcolor=magenta]{hyperref}
\usepackage{cleveref}

\DeclareMathOperator{\TCL}{TCL}

\makeatletter
 
 \@addtoreset{equation}{section}
\makeatother

\newtheorem{theorem}{Theorem}

\newtheorem{proposition}[theorem]{Proposition}

\theoremstyle{remark}
\newtheorem{remark}[theorem]{Remark}

\numberwithin{equation}{section}

\usepackage{lineno}

\title{The total chord length of maximal outerplanar graphs}
\author{Haley Broadus and Elena Pavelescu}

\begin{document}

\begin{abstract}
We consider embeddings of maximal outerplanar graphs whose vertices all lie on a cycle $\mathcal{C}$ bounding a face. Each edge of the graph that is not in $\mathcal{C}$, a chord,  is assigned a length equal to the length of the shortest path in $\mathcal{C}$ between its endpoints. We define the total chord length of a graph as the sum of lengths of all its chords. For each order $n\ge 5$, we find outerplanar graphs whose total chord length is the minimum among all graphs of the same order, and graphs whose total chord length is the maximum among all graphs of the same order. We give a complete characterization of the graphs that attain the maximum total chord length. We show that every integer value in the interval determined by the minimum and maximum values is the total chord length of a maximal outerplanar graph of the same order.
\end{abstract}
\maketitle

\section{Introduction}
An \emph{outerplanar} graph is a graph that has a planar embedding such that all the vertices belong to the same face of the embedding. A \emph{maximal outerplanar} graph is an outerplanar graph that fails to remain outerplanar after the addition of any one edge. 
Outerplanarity is a hereditary property
and it is characterized by the forbidden minors $K_4$ and $K_{3,2}$ \cite{CH}.
The \emph{order} of a graph represents its number of vertices.  
A maximal outerplanar graph of order $n$ can be drawn as an $n-$cycle $\mathcal{C}$ together with a triangulation of the bounded plane region determined by $\mathcal{C}$.
In this article, a maximal outerplanar graph means a drawing which includes $\mathcal{C}$ and a triangulation of the bounded region.
The edges of the graph that do not belong to the cycle $C$ are called \textit{chords}. 
 A maximal outerplanar graph of order $n$ has $n-3$ chords and each chord is assigned a \emph{length} equal to the shortest path between its endpoints included in the cycle $\mathcal{C}$. A chord can have length as small as 2 and as large as $\lfloor \frac{n}{2}\rfloor$. A chord of length 2 is called an \textit{ear} of the graph.  A chord of length $\frac{n}{2}$ is called a \textit{diameter} of the graph.
The \textit{total chord length (TCL)} of a maximal outerplanar graph equals the sum of the lengths of all its chords. 
\begin{figure}[ht]
\centering
\includegraphics[width=0.4\textwidth]{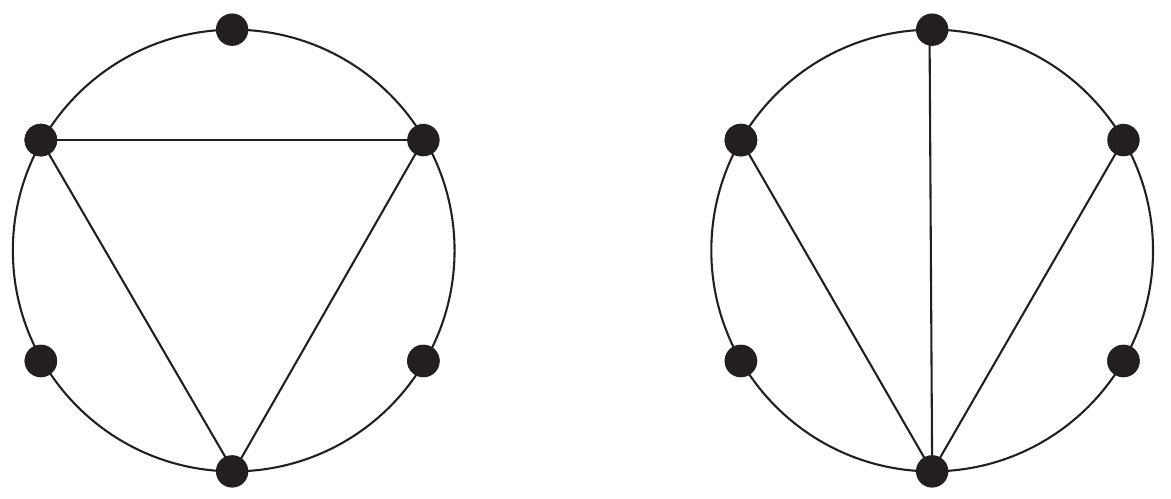}
\caption{\label{fig:example6} The two maximal outerplanar graphs on six vertices. The graph on the left has three chords of length 2 and $\TCL=6$, and the graph on the right has two chords of length 2, one chord of length 3, and  $\TCL=7$.  }
\end{figure}

In this article, we derive explicit formulas for the maximum and minimum total chord lengths over all maximal outerplanar graphs of a given order.  These values are given in Theorems \ref{thm:minimum} and \ref{thm:maximum}. We prove Theorems \ref{thm:minimum} and \ref{thm:maximum} in Sections \ref{secmin}   and \ref{secmax}, respectively, and in each case we give specific examples of graphs for which these extreme values are realized. 
\begin{theorem} 
Let $k\ge 0$. For $3\cdot 2^k \le n \le 3\cdot 2^{k+1}$, the minimum $\TCL$ 
among all maximal outerplaner graphs of order $n$ is: 
\[ 
min\TCL(n)=  n(k+2)-3\cdot 2^{k+1}.
\]
\label{thm:minimum}
    \end{theorem}
    
\begin{theorem} 
Let $n\ge 5$. The maximum $\TCL$ among all maximal outerplanar graphs of order $n$  is 
\[ 
max\TCL(n)= \left\{
\begin{array}{ll}
      \frac{n^2-9}{4},& \textrm{if } n \textrm{ odd} \\
 \frac{n^2-8}{4},& \textrm{if } n \textrm{ even.} \\\end{array} 
\right. 
\]
\label{thm:maximum}
    \end{theorem}
\noindent In Section \ref{secmax}, we also give a complete
characterization of the graphs that attain the maximum total chord length:  graphs which have exactly two ears. In Section \ref{minmax}, we show that for every order $n\ge 5$ and every $min\TCL(n)\le l\le max\TCL(n)$, there exists a maximal outerplanar graph of order $n$ whose total chord length equals $l$. 

The problem of triangulating a set of points has been studied extensively in computational geometry. The total chord length  provides a new combinatorial parameter that reflects the structural properties of outerplanar graph triangulations. The problem of finding the minimum TCL is a special case of the minimum-weight triangulation (MWT) problem first posed in \cite{DG}.
The input to the MWT  problem is a set of points $S$ in the plane, and the desired output is a triangulation of  $S$ (a maximal planar straight-line graph with vertex set $S$)  that minimizes the total Euclidean edge length \cite{BE}.
In \cite{MR}, Mulzer and Rote prove that  the MWT problem is NP-hard.
A triangulation of minimum total chord length for a \textit{regular $n$-gon} is also a minimal-weight triangulation for the set of its vertices.
This is because the edges of the outer circle are included in the triangulation \cite{BYZ}, and the larger the chord length of a chord, the larger its Euclidean length. 
In this special case, the problem can be solved in $\mathcal{O}(n^3)$ time using a dynamic programming algorithm by Gilbert \cite{G}, \cite{CLRS}.

\section{Background and Notation}

\begin{figure}[ht]
\centering
\includegraphics[width=0.55\textwidth]{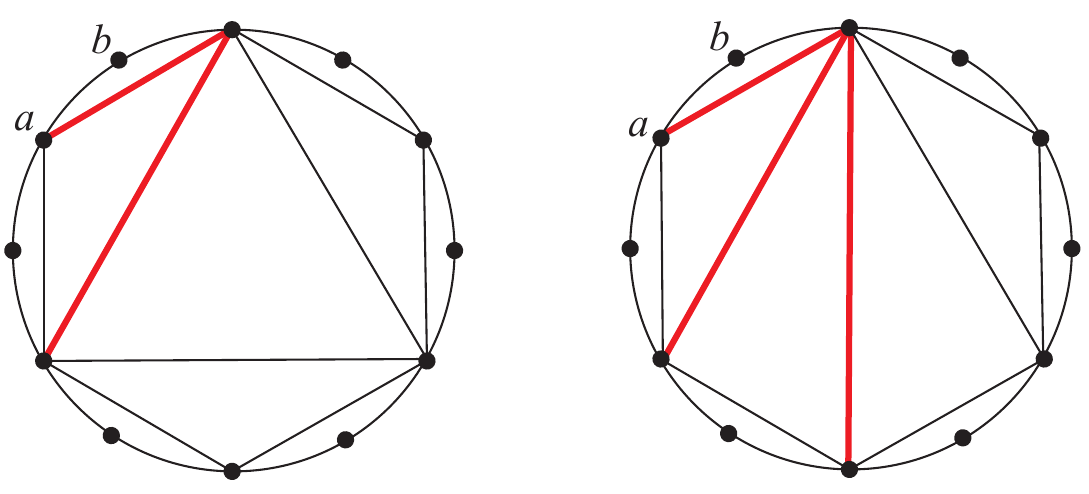}
\caption{\label{fig:example12} Two maximal outerplanar graphs with twelve vertices. }
\end{figure}

In a maximal outerplanar graph, each ear corresponds to a vertex of degree 2 of the graph. 
The graphs in Figure \ref{fig:example6} have three and two ears, respectively. 
The following theorem, frequently attributed to  Meisters \cite{M} is called the ``two ears theorem" and can be proved by induction on the number of vertices.
\begin{theorem} A maximal outerplanar graph with $n\ge 5$ vertices has at least two ears.  
\label{thm:2ears}
\end{theorem}
Consider an outerplanar graph $G$ of order $n$ and let $\mathcal{C}$ be its $n-$cycle.
For an edge $ab \in E(\mathcal{C})$, we define $n_G(ab)$ as the number of chords that ``layer"  $ab$. More precisely,
\begin{itemize}
    \item If $G$ does not contain a diameter:  $n_G(ab) := $ number of chords in $G$ whose length is computed along a path in $\mathcal{C}$ that includes $ab$.
    \item If $G$ does contain a diameter: $n_G(ab) := $ $ \frac{1}{2}$  + number of chords in $G$ whose length is computed along a path in $\mathcal{C}$ that includes $ab$, except the diameter.
\end{itemize}
In Figure \ref{fig:example12}, $n_G(ab)=2$ for the graph on the left, and $n_G(ab)=2\frac{1}{2}$ for the graph on the right. We note that $$\TCL(G)= \sum_{ab\in E(\mathcal{C})} n_G(ab).$$

We also note that for a graph $G$ of order $n$ and an edge $ab\in E(\mathcal{C}),$
\[ n_G(ab) \le  \frac{n-3}{2}. \]
This is because when $n$ is odd, at most half of the $n-3$ chords have their length computed along a path in  $\mathcal{C}$ that includes $ab$. And when $n$ is even, at most half of $n-4$ chords plus a diameter have their length computed along a path in  $\mathcal{C}$ that includes $ab$.
An \textit{edge contraction} means that the endpoints of the edge are identified to one vertex that inherits all the neighbors of the two endpoints. Loops and double edges created by the identification of the two endpoints are deleted.
The contraction of a path within a graph means contracting each edge within that path.
 \textit{Subdividing the edge} $ab$ means deleting the edge $ab$ and adding a new vertex $c$ together with edges $ac$ and $bc$.
 \noindent By Theorem \ref{thm:2ears}, each graph of order $n \ge 5$ has at least two ears. Consequently, for $n \ge 5$, each graph of order $n + 1$ can be obtained from a graph of order $n$ by subdiving an edge $ab\in E(\mathcal{C})$ and adding the chord $ab$.

 \begin{remark}
Assume $G'$ is a graph obtained from $G$ by subdividing an edge $ab \in E(\mathcal{C})$ in $G$ and adding the chord $ab$.\\
Then 
\[
\TCL(G') = 
\begin{cases} 
\TCL(G) + 2 + n_G (ab), & \text{if } G \text{ has no diameter}, \\
\TCL(G) + \frac{3}{2} + n_G (ab), & \text{if } G \text{ has a diameter}.
\end{cases}
\]
\label{remark:subdivide}
\end{remark}
\begin{proof}
If $G$ has no diameter, then $\TCL$ increases by 2 due to the new chord $ab$, and it increases by 1 for each chord of $G$ that layers $ab$. In total, there is an increase of $2+n_G(ab)$. If $G$ has a diameter, then $\TCL$ increases by 2 due to the new chord $ab$, and it increases by 1 for each chord of $G$ that layers $ab$ and is not the diameter. In total, there is an increase of $2+(n_G(ab)-\frac{1}{2})= \frac{3}{2}+n_G(ab)$. The $-\frac{1}{2}$ comes from the loss of the diameter that layered $ab$ in $G$.
\end{proof}
    


\section{Minimum Total Chord Length}
\label{secmin}
This section starts with an algorithm to construct a maximal outerplanar graph for each order $n\ge 4$. Then the TCL for the graph obtained through this algorithm is computed and shown to be the smallest possible for all graphs of the same order. 

\begin{figure}[ht]
\centering
\includegraphics[width=0.55\textwidth]{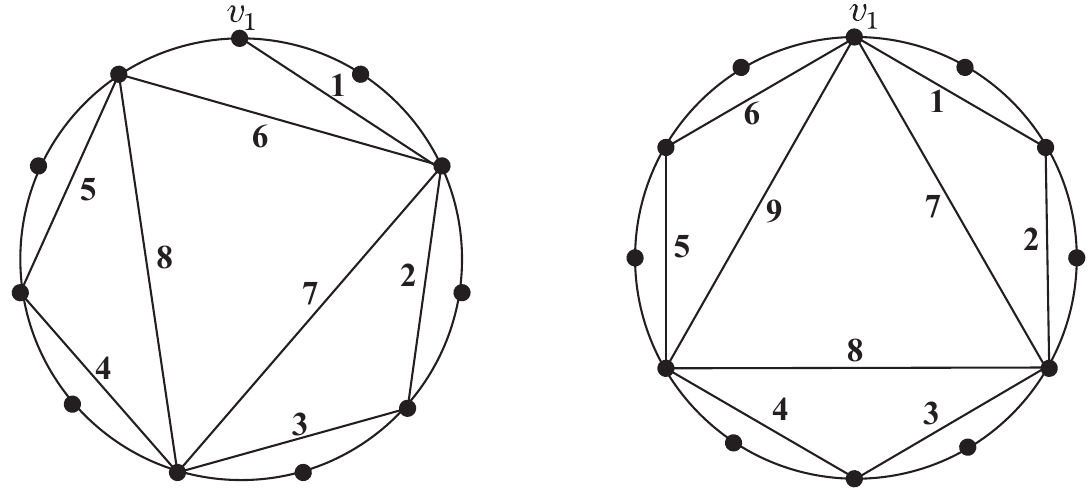}
\caption{\label{fig:algorithm} Maximal outerplanar graphs of order eleven (left) and twelve (right). }
\end{figure}

For $n\ge 4$, consider an $n-$cycle $\mathcal{C}$ whose vertices are labeled $v_1$, $v_2, \ldots, v_n$, in clockwise order, as they appear in $\mathcal{C}$. Construct a maximal outerplanar graph as follows, where edges are added only within the bounded region determined by $\mathcal{C}$: (1) Include the shortest chord possible with starting endpoint at $v_1$  and final endpoint situated after $v_1$ in the clockwise direction around $\mathcal{C}$. (2) Make the final endpoint of the previously added chord the starting endpoint of a shortest chord whose final endpoint is positioned in the clockwise direction from its starting endpoint. (3) Continue to add chords of minimum length this way until $n-3$ chords have been added.  Two graphs constructed following this algorithm are shown in Figure \ref{fig:algorithm}. The chords are labeled according to the order in which they were added. 

We will show that for each order, this algorithm gives a maximal outerplanar graph that minimizes the total chord length. We call graphs constructed in this manner \emph{greedy graphs}. Note that the greedy graph triangulation is what is known as the balanced hierarchical decomposition or Dobkin–Kirkpatrick hierarchy \cite{DK}.

\begin{proposition}
Let $k\ge 0$. For $3\cdot 2^k \le n \le 3\cdot 2^{k+1}$, the  total chord length of the maximal outerplanar graph $G_n$ of order $n$ obtained with the algorithm above is 
\[ 
\TCL(G_n)=  n(k+2)-3\cdot 2^{k+1}. \]
\label{thm:min}
    \end{proposition}
    \begin{proof}
    For $k= 0$, $3\le n\le 6$, we check that $\TCL(G_n)=2n-6,$ which matches the statement.
    Let $k\ge 1$. For $n = 3\cdot 2^k$, the graph $G_n$ contains a triangle whose three edges are chords of length $2^k$. For each $ab\in E(\mathcal{C})$,  $n_G (ab)=k$, and 
\[\TCL(G_n) = nk = n (k+2) -2n = n(k+2)-3\cdot 2^{k+1} .\] 
For each $3\cdot 2^k \le  n \le  3\cdot 2^{k+1}$, a greedy graph of order $n$ is obtained from the greedy graph on $3\cdot 2^k$ vertices by subdividing edges of $\mathcal{C}$ in counterclockwise order as they appear in $\mathcal{C}$.
Figure \ref{fig:greedy1216} contains greedy graphs of orders  12 through 16 ($k=2$). Graphs up to order 24 can be obtained by subdividing one by one each of the remaining eight edges of the original 12-cycle. The label $v_1$ marks the vertex at the beginning of the algorithm that gives the greedy graph of that order. 
If $G'$ is obtained from $G$ by one such edge subdivision, then \[\TCL(G')=\TCL(G)+(k+2).\]

\begin{figure}[ht]
\centering
\includegraphics[width=1\textwidth]{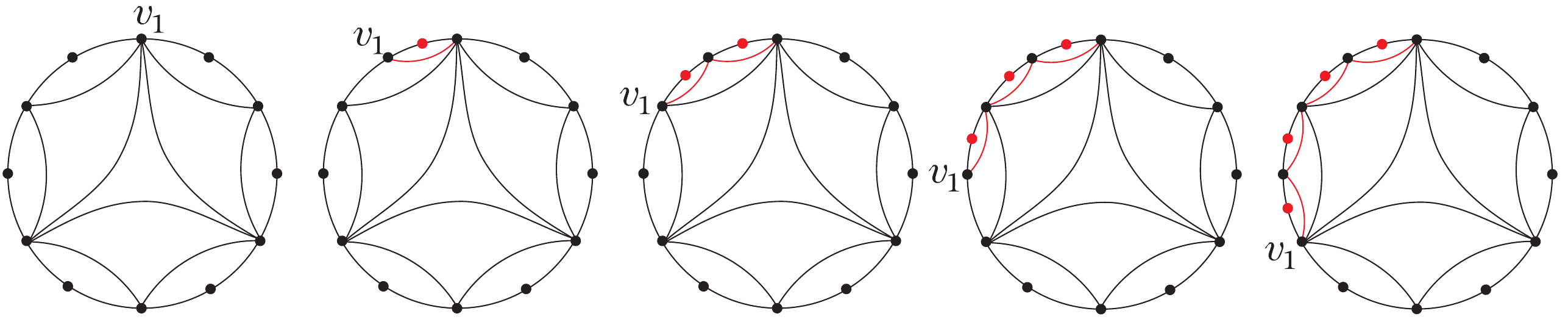}
\caption{\label{fig:greedy1216} Maximal outerplanar graphs of orders twelve through sixteen. }
\end{figure}

\noindent This is because a chord of length 2 is added, and the length of each of $k$ chords which layered the subdivided edge is increased by 1. 
In turn, for $3\cdot 2^k\le n \le 3\cdot 2^{k+1}$,
\[
\begin{aligned}
    \TCL(G_n) &= \TCL(G_{3\cdot2^k})+ (n-3\cdot 2^k)(k+2) \\
              &= 3\cdot 2^k \cdot k  + (n-3\cdot 2^k)(k+2) \\
              &= n(k+2) -3\cdot 2^{k+1}.
\end{aligned}
\]

    \end{proof}

  Next, we show that for each order $n\ge 4$, the algorithm gives a maximal outerplanar graph that minimizes the total chord length.


\begin{theorem}
    For a fixed integer $ k \geq 1$, let $\theta(k)$ represent the largest value for which there exists a maximal outerplanar graph of order $\theta(k)$ such that  $n_G (ab) \leq k$ for all $ab \in E(\mathcal{C})$. Then $\theta(k) = 3 \cdot 2^k.$
    \label{thm:theta}
\end{theorem}
\begin{proof}
 
\noindent For a maximal outerplanar graph, consider the maximum among the lengths of all its chords. No edge of $E(\mathcal{C})$ contributes to the length of more than one chord of maximal length. Otherwise, the chords of maximal length will be ``nested", contradicting maximality. This means all chords of maximal length are adjacent to a common triangle within the triangulation of the disk, and therefore, there are at most three such chords. The subset of edges of $E(\mathcal{C})$ layered by one chord of maximal length is disjoint from the subset of edges of $E(\mathcal{C})$ layered by a different chord of maximal length. This implies a graph of the largest order satisfying $n_G(ab) \le  k$ for all $ab \in E(\mathcal{C})$ has three chords of the greatest length.
   
    We proceed by induction on $k$. The base case, $k=1$, is true since the largest order $n$ for which $n_G(ab) = 1$ is $n = 6$. 
    Assume the statement holds for $k$ and let $G$
be a graph  of order $n = 3 \cdot 2^k$ with $n_G(ab) \leq k$ for each $ab \in E(\mathcal{C})$. 
    For each edge $ab \in E(\mathcal{C})$, subdivide $ab$ and add a chord from $a$ to $b$, creating a graph $G'$ of order $3 \cdot 2^{k+1}$. The construction ensures $n_{G'}(ac) \leq k+1$ for each $ac \in E(\mathcal{C}')$, implying $\theta(k+1) \geq 3 \cdot 2^{k+1}$.

\begin{figure}[ht]
\centering
\includegraphics[width=0.25\textwidth]{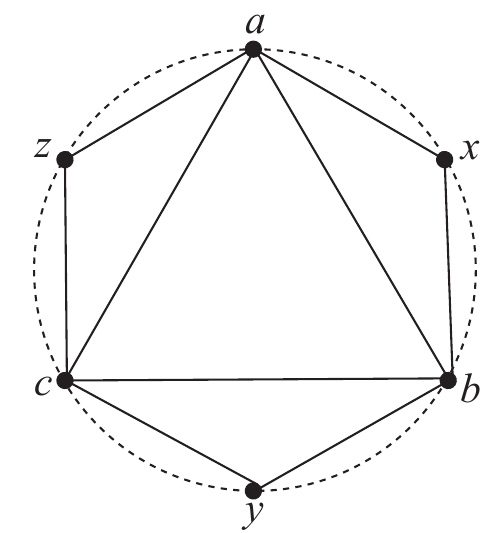}
\caption{\label{fig:triangle} The graph $H$ is  obtained by deleting the triangle $abc$ and contracting the three shortest arcs among $ax$, $bx$, $by$, $cy$, $az$, and $cz$. }
\end{figure}

    Assume $\theta(k+1) > 3 \cdot 2^{k + 1}$. The graph of greatest order for which $n_G(ab) \le k+1$ for all $ab\in E(\mathcal{C})$ has three chords of largest length. Call the three endpoints of these three chords $a,b,$ and $c$ and let $T$ be the triangle they determine. The chords $ab$, $ac$, and $bc$ each belong to a triangle different from $T$. That is, there exist vertices $x, y, z \in V(G)$ such that $abx$, $acy$, and $bcz$ are also triangles of the graph. See Figure \ref{fig:triangle}. The vertices $a, b, c, x, y,$ and $z$ determine six arcs of $\mathcal{C}$. We construct a graph $H$ by deleting the triangle $T$ and contracting to one vertex each of the three shortest arcs among $ax$, $bx$, $by$, $cy$, $az$, and $cz$.
The graph $H$ is of order greater than $3 \cdot 2^{k+1} \cdot \frac{1}{2} = 3 \cdot 2^k$ and $n_H(ab) \leq k$. This contradicts the induction hypothesis. This means $\theta(k+1) = 3 \cdot 2^{k+1}.$
\end{proof}

\begin{theorem}
    For $n\ge 3$, the greedy graph $G_n$ has minimum total chord length among all maximal outerplanar graphs of order $n$.
    \end{theorem}
    
    \begin{proof}
   For orders 3, 4, and 5, there is a unique maximal outerplanar graph; thus, the statement holds. The graph $G_6$ has the minimum chord length among all maximal planar graphs with six vertices. Assume by induction that the greedy graph $G_n$ has minimum total chord length among all graphs of order $n$.
        For $n \geq 6$, and $3 \cdot 2^k \leq n < 3 \cdot 2^{k+1}$,
        \[\TCL(G_{n+1}) = \TCL(G_n) + k + 2.\]
By Theorem \ref{thm:theta}, for any other graph $G'$ of order $n+1$, where $3 \cdot 2^k \leq n < 3 \cdot 2^{k+1}$, there exists an edge $ab$ such that $n_{G'}(ab)\ge k+1$. Without loss of generality, assume the graph $G'$ is obtained from a graph $G$ of order $n$ by subdividing an edge $ac$, adding the new vertex $b$, and adding the chord $ac$. We explain below that  
 \[\TCL(G') \geq \TCL(G) + k + 2,\] which implies 
        \[\TCL(G') \geq  \TCL(G) + k + 2 \ge \TCL(G_n)+k+2 = \TCL(G_{n+1}),\]
        and thus for each $n\ge 6$, the greedy graph $G_n$ has the minimum TCL among all graphs of order $n$.

        If the graph $G$ has a diameter, then  $k+1\le n_{G'}(ab)=n_G(ac)+\frac{1}{2},$ thus $n_G(ac)\ge k+\frac{1}{2}.$
        If the graph $G$ doesn't have a diameter, and the graph $G'$ has a diameter, then $k+1\le n_{G'}(ab)=n_G(ac)+ \frac{1}{2},$ thus $n_G(ac)\ge k+\frac{1}{2}.$
        Lastly, if neither graphs $G$ and $G'$ have a diameter, then $k+1\le n_{G'}(ab)=n_G(ac)+ 1,$ thus $n_G(ac)\ge k.$ In all cases, $n_G(ac)\ge k.$ By Remark \ref{remark:subdivide}, since TCL$(G')$ is an integer,
        \[\TCL(G') \geq \TCL(G) + k + 2.\]
    \end{proof}
\begin{remark} The process of adding a chord involves moving to the next available vertex in clockwise order. A maximal outerplanar graph with \(n\) vertices has exactly \(n-3\) chords. Therefore, the overall time complexity of the algorithm is \(O(n)\).
\end{remark}
\begin{remark} In the algorithm outlined in this section, a chord of minimum length is added at each step according to a specified method. However, if chords of minimum length are added without constraining their endpoints, the resulting graph may not necessarily have the minimum total chord length. It might be valuable to find other algorithms yielding graphs with the minimum total chord length.
    \end{remark}
    

\section{Maximum Total Chord Length}
\label{secmax}
In this section we describe a particular maximal outerplanar graph, and then we compute its total chord length. We show that this TCL is the largest possible for all graphs of the same order. 

For $n\ge 5$, consider the graph $F_n$ with $n$ vertices obtained from the cycle of length $n$ whose vertices are labeled $v_1, v_2, \ldots, v_n$ in clockwise direction by adding the chords $v_1v_3, v_1v_4, \ldots, v_1v_{n-1}.$ 
The resulting graph resembles a hand fan, and is referred to as a \emph{fan graph}.

\begin{figure}[ht]
\centering
\includegraphics[width=0.6\textwidth]{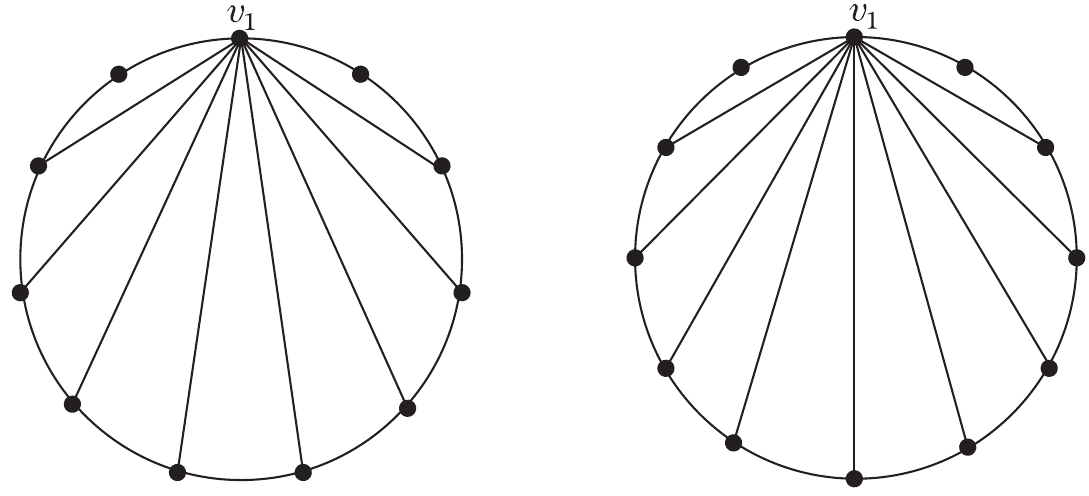}
\caption{\label{fig:greedy} Fan graphs of orders eleven(left)  and twelve(right). }
\end{figure}

\begin{proposition}
Let $n\ge 5.$ The fan graph $F_n$ with $n$ vertices has total chord length \[ 
\TCL(F_n)= \left\{
\begin{array}{ll}
      \frac{n^2-9}{4},& n \textrm{ odd} \\
 \frac{n^2-8}{4},& n \textrm{ even} \\\end{array} 
\right. 
\]
\end{proposition}
\begin{proof}
For $n$ odd, the graph has two chords of length $i$ for each $i=2, 3, \ldots, \frac{n-1}{2}$. In this case 
\[ \TCL(F_n) = 2 \cdot \sum_{i=2}^{(n-1)/2} i  = 2\cdot \big(\frac{\frac{n-1}{2}(\frac{n-1}{2}+1)}{2} -1\big)= 2\cdot \big(\frac{\frac{n-1}{2}\cdot \frac{n+1}{2}}{2} -1\big)= \frac{n^2-1}{4}-2=\frac{n^2-9}{4}.\]

For $n$ even, the graph has a diameter and two chords of length $i$ for each $i=2, 3, \ldots, \frac{n-2}{2}$. In this case 
\[ \TCL(F_n) = \frac{n}{2}+ 2 \cdot \sum_{i=2}^{(n-2)/2} i = \frac{n}{2} + 2\cdot \big(\frac{\frac{n-2}{2}(\frac{n-2}{2}+1)}{2} -1\big)= \frac{n}{2}+ \frac{n-2}{2}\cdot\frac{n}{2}-2= \frac{n^2-8}{4}.\]
\end{proof}
Notice that there is inherent symmetry within the fan graphs: for every chord of length \(l\) on the left-hand side, there is a corresponding chord of length \(l\) on the right-hand side. Additionally, observe that the chord length grows in increments of 1, starting at length 2.

\begin{proposition}\label{prop:fan}
The fan graph $F_n$ has the maximum total chord length among all maximal outerplanar graphs with $n$ vertices.
\end{proposition}
\begin{proof}
\noindent 
Assume the graph $G'$ with $n+1$ vertices is obtained from the graph $G$ of order $n$ by subdividing an edge $ab \in E(\mathcal{C})$ and adding the chord $ab$. Then by Remark \ref{remark:subdivide},
\[
\TCL(G') = 
\begin{cases} 
\TCL(G) + 2 + n_G(ab), & \text{if } G \text{ has no diameter}, \\
\TCL(G) + \frac{3}{2} + n_G(ab), & \text{if } G \text{ has a diameter}.
\end{cases}
\]
Since $n_G(ab) \le \frac{n-3}{2},$ we have that
\begin{equation}
\TCL(G') \le  
\begin{cases} 
\TCL(G) + \frac{n+1}{2}, & \text{if } G \text{ has no diameter}, \\
\TCL(G) + \frac{n}{2}, & \text{if } G \text{ has a diameter}.
\end{cases}
\label{EquationTLC1}
\end{equation}

We also have
\begin{equation}
\TCL(F_{n+1})=  
\begin{cases} 
\TCL(F_n) + \frac{n+1}{2}, & \text{if } n \text{ odd}, \\
\TCL(F_n) + \frac{n}{2}, & \text{if } n \text{ even}.
\end{cases}
\label{EquationTLC2}
\end{equation}
The fan graph $F_5$ has the maximum total chord length, as the only graph of order 5. Assume by induction that $F_n$ has maximum total chord length among graphs of order $n.$ In all cases (whether $n$ is odd or even and whether $G$ has a diameter or not), equations \ref{EquationTLC1} and \ref{EquationTLC2} imply  $\TCL(G') \le \TCL(F_{n+1})$, and the proof is complete by induction.
\end{proof}
    
\begin{theorem}\label{thm:max2ears}
A maximal outerplanar graph has maximum total chord length if and only if it has exactly two ears. 
\end{theorem}
\begin{proof}
To prove the direct implication, assume $G$ has order $n$ and at least three ears. 
The graph $G$ has one, two, or three chords of maximum length. 
In each case, $G$ has a non-chord $ab\in E(\mathcal{C})$ such that $n_G(ab)$ is less than the maximum possible value $\frac{n-3}{2}$. See Figure \ref{fig:3ears}. When there are at most two longest chords, the longest chord(s) create two sides of the graph. One of these sides will contain two ears. A cycle-edge $ab$ layered by one of the two ears is not layered by the other ear, and thus $n_G(ab) <\frac{n-3}{2}$.
When there are three longest chords, the graph is evenly divided into three sections, with each section containing at least one ear. For a cycle-edge $ab$ layered by an ear, $n_G(ab) \leq \frac{n-3}{3} < \frac{n-3}{2}.$
This means that $G$ is obtained from a graph $\overline{G}$ of order $n-1$ through a subdivision of a non-chord, and $\TCL(G)-\TCL(\overline{G})<\TCL(F_n)-\TCL(F_{n-1}).$ This implies $G$ doesn't have maximum total chord length. 

\begin{figure}[ht]
\centering
\includegraphics[width=0.9\textwidth]{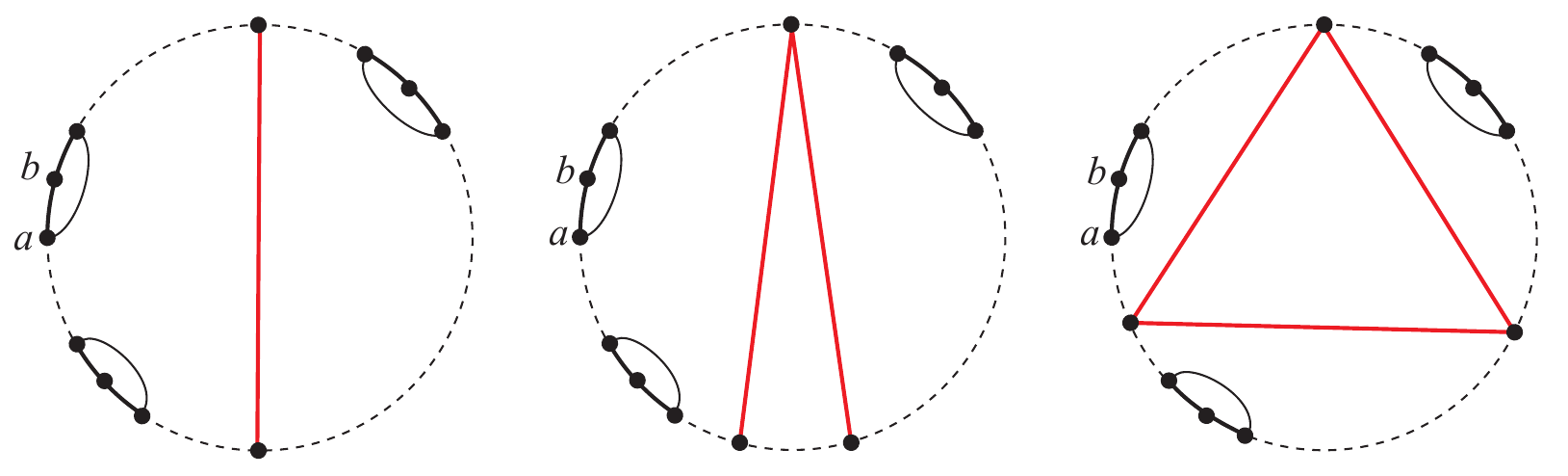}
\caption{\label{fig:3ears} Maximal outerplanar graphs with one, two, and three chords of maximum length and non-chord $ab$ with less than maximum $n_G.$}
\end{figure}

To prove the reverse, consider a graph $G$ that has exactly two ears. 
A chord of maximum length layers at least one ear.
Since the graph has exactly two ears, then $G$ has at most two chords of maximum length. Assume $G$ has one chord $d$ of maximum length. 
This means $n$ is even. 
The chord $d$, together with each half of the cycle $\mathcal{C}$ that the endpoints of $d$ determine are boundary cycles for two maximal outerplanar graphs $H_1$ and $H'_1$ of order $\frac{n}{2}+1$.
These graphs are highlighted in red ($H_1$) and in blue ($H_1'$) in Figure  \ref{fig:maximal}.
As maximal outerplanar graphs, each of these graphs have two ears. Since only one of these two ears is also an ear of $G$, the second ear of $H_1$ layers  $d$ and has length $\frac{n}{2}-1$ in $G$.  Call this chord $d_1.$ See Figure \ref{fig:maximal} (left).
The chord $d_1$, together with the path along $\mathcal{C}$ which gives its length in $G$ form the boundary cycle for an outerplanar graph $H_2$ of order $\frac{n}{2}.$ The graph $H_2$ has two ears, one of which layers $d_1$  and has length  $\frac{n}{2}-2$ in $G$. Call this chord $d_2.$ 
Continue until a chord of length 3 is constructed. 
This reasoning can be replicated for the graph $H_1'$.
The graph $G$ has one chord of length $\frac{n}{2},$ and two chords each of lengths $2, 3, \ldots,  \frac{n}{2}-1.$ Its total chord length is equal to that of the fan graph.

If the graph has two chords of maximum length $d$ and $d'$ ($n$ odd) the same reasoning applies. See Figure \ref{fig:maximal} (right).

\begin{figure}[ht]
\centering
\includegraphics[width=0.6\textwidth]{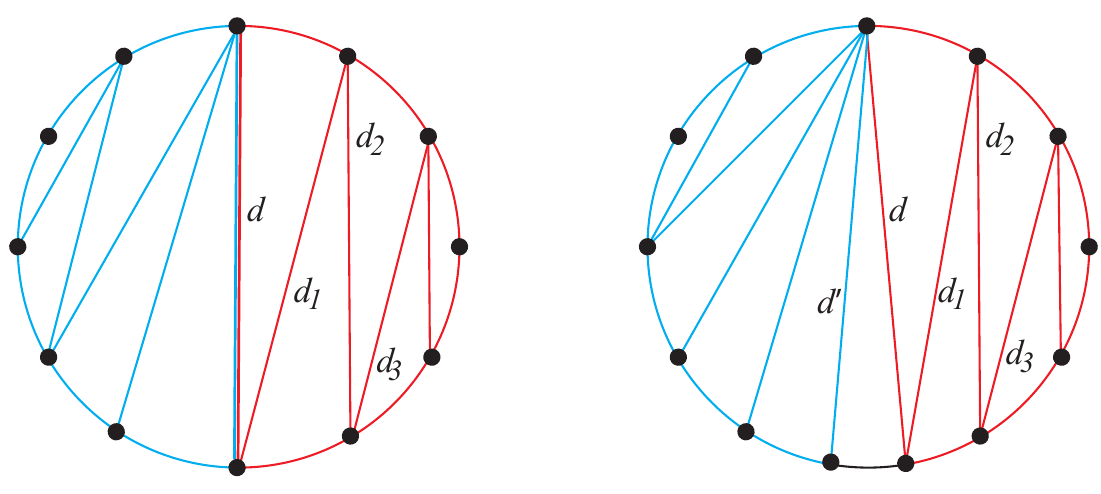}
\caption{\label{fig:maximal} Maximal outerplanar graphs of even order (left) and odd order (right) with maximum total chord length. }
\end{figure}
 
\end{proof}

\section{All possible values for the total chord length}
\label{minmax}
 Let $G_n$ represent the greedy graph on $n$ vertices and $F_n$ represent the fan graph on $n$ vertices as described in Sections \ref{secmin} and \ref{secmax}.
\begin{theorem} For every $n\ge 5$ and $\TCL(G_n)\le l\le \TCL(F_n)$, there exists a maximal outerplanar graph of order $n$ with total chord length equal to $l$.
\label{thm:flip}
\end{theorem}

\begin{proof}
 We begin with the fan graph $F_n$ and describe an algorithm that transitions to the greedy graph $G_n$ by a sequence of flips, where each step permanently adds chords of decreasing length and shifts the anchoring vertex. Between iterations, additional flips can be made to control the change in total chord length more finely. Specifically, for $n = 2^k \cdot m$ ($m$ odd), at the $r$-th iteration, certain flips decrease the total chord length by up to $2^{r-1}$, possibly leaving gaps. We show how to fill these gaps using additional flips that increase the total chord length in increments from 1 to $2^{r-1}-1$, thereby realizing all values between $\TCL(G_n)$ and $\TCL(F_n)$.
 
\begin{figure}[ht]
\centering
\includegraphics[width=0.9\textwidth]{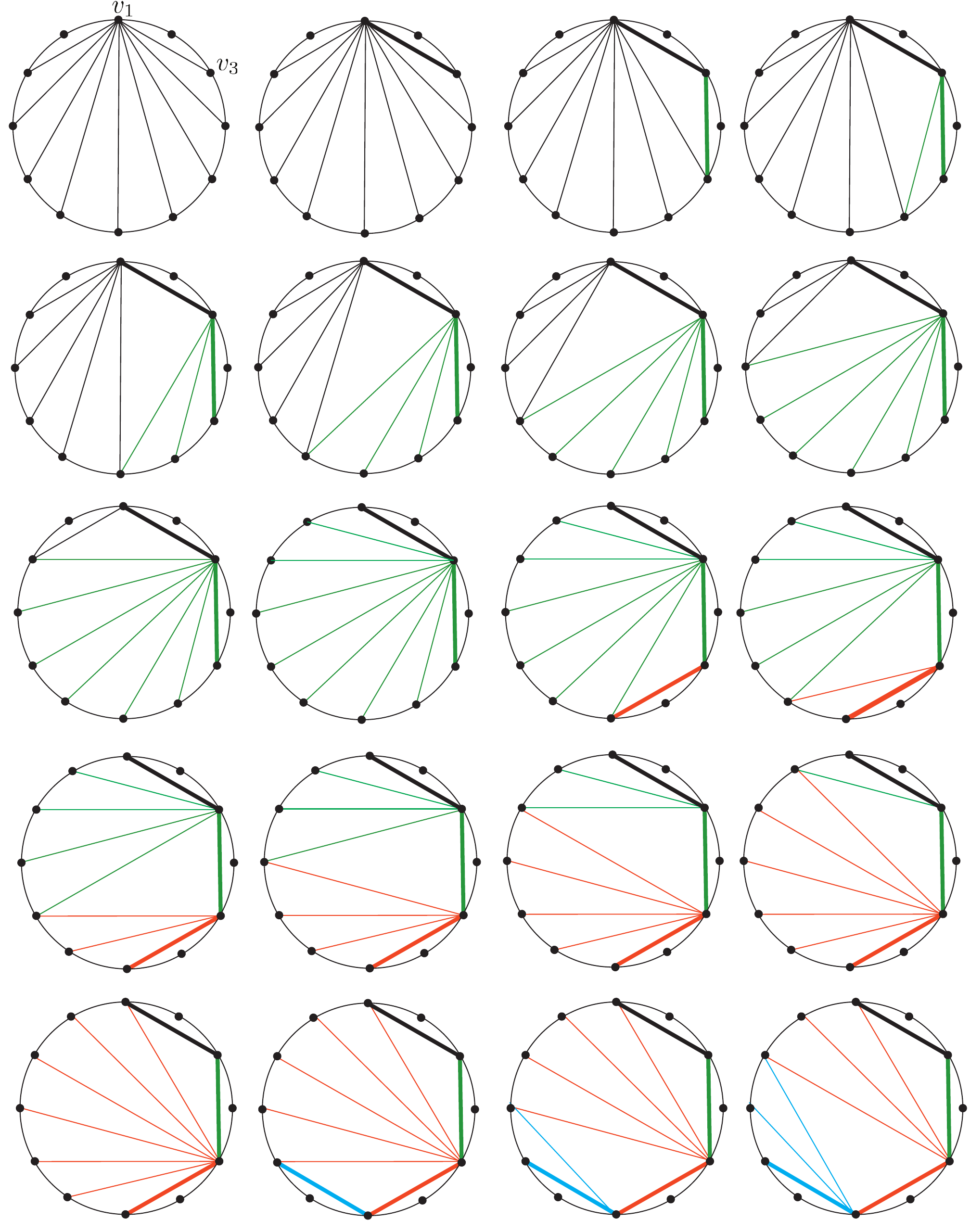}
\caption{\label{fig:algo1} Re-anchoring algorithm applied to the fan graph on twelve vertices.}
\end{figure}
\noindent The algorithm for order 12 is presented in Figures \ref{fig:algo1} and \ref{fig:algo2}. The bold edges in each graph represent edges that are part of the greedy graph. Once these ``permanent" edges are formed, they persist through to the final graph.
The algorithm can be summarized as follows:
 \begin{figure}[ht]
\centering
\includegraphics[width=0.9\textwidth]{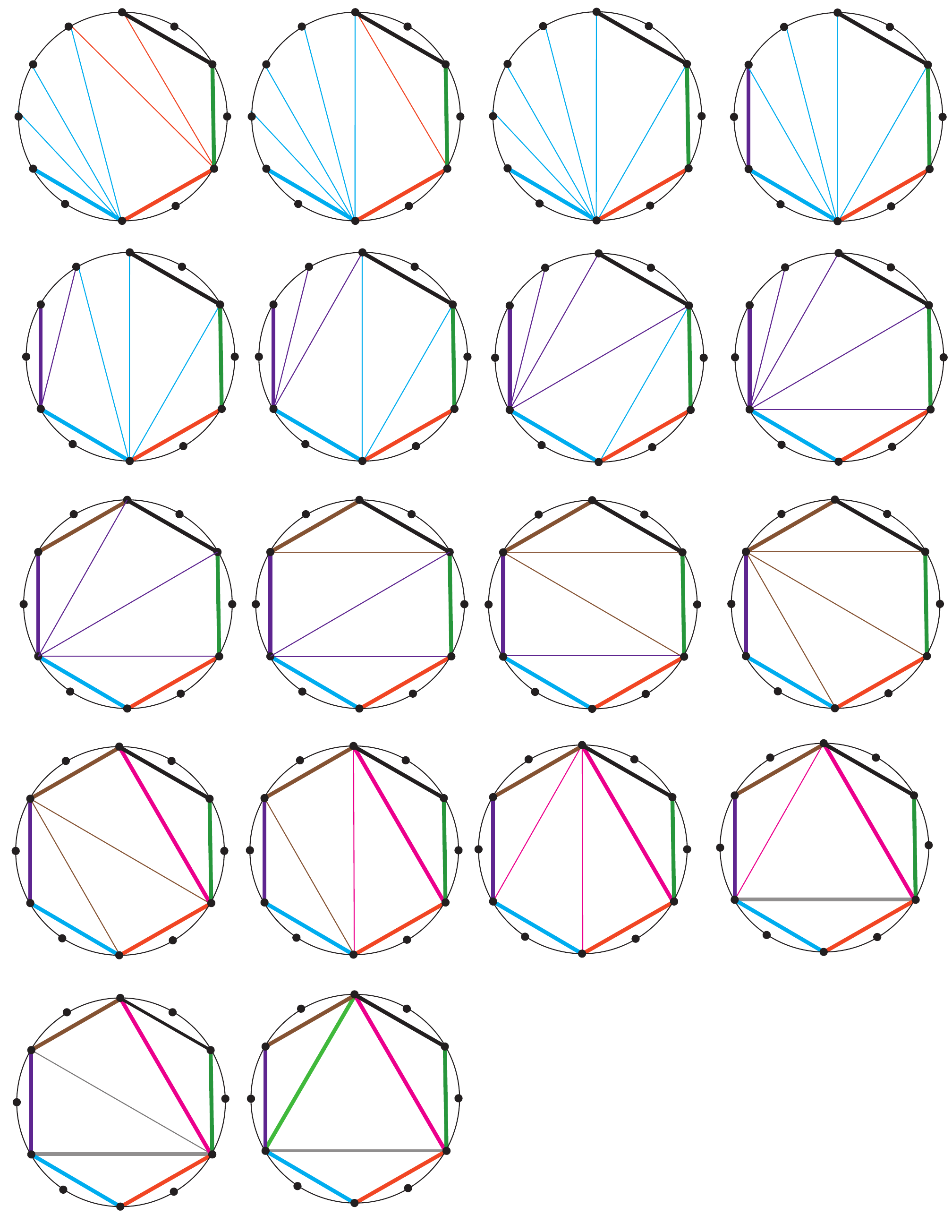}
\caption{\label{fig:algo2} Re-anchoring algorithm applied to the fan graph on twelve vertices (continued)}
\end{figure}
Start with the graph \( F_n \) whose chords are all incident to $v_1$. Mark as permanent the chord $v_1v_3$. One by one, replace the other chords incident to $v_1$ by chords incident to $v_3$, as in the first three rows of Figure \ref{fig:algo1}. In this process, \( v_3 \) serves as a new ``anchor," and the chord $v_3v_5$ is a permanent chord.
Next, replace the chords incident to $v_3$ (except the two permanent chords) one by one with chords incident to $v_5$. The permanent chord  \( v_5v_7 \) is added, and \( v_5 \) becomes the new anchor. 
Since a permanent chord is added each time the anchoring vertex is changed, with each change in anchor, one vertex is eliminated from the set of possible future anchors, and each subsequent anchor necessitates one less chord. 
This process will end, and by way of construction, the final graph is the greedy graph. 


It remains to prove that for each value  $\TCL(G_n)\le l\le \TCL(F_n)$, there is a graph of order $n$ whose total chord length is $l$.
Not all values $l$ are achieved by the graphs given by the algorithm. We will show how these missing values can be obtained by modifying the graphs produced by the algorithm.
For $r\ge 1$, we say the $r$th iteration of the algorithm is completed when all chords of length $2^{r+1}-1$ or smaller have been added as permanent chords. 

\begin{figure}[ht]
\centering
\includegraphics[width=0.9\textwidth]{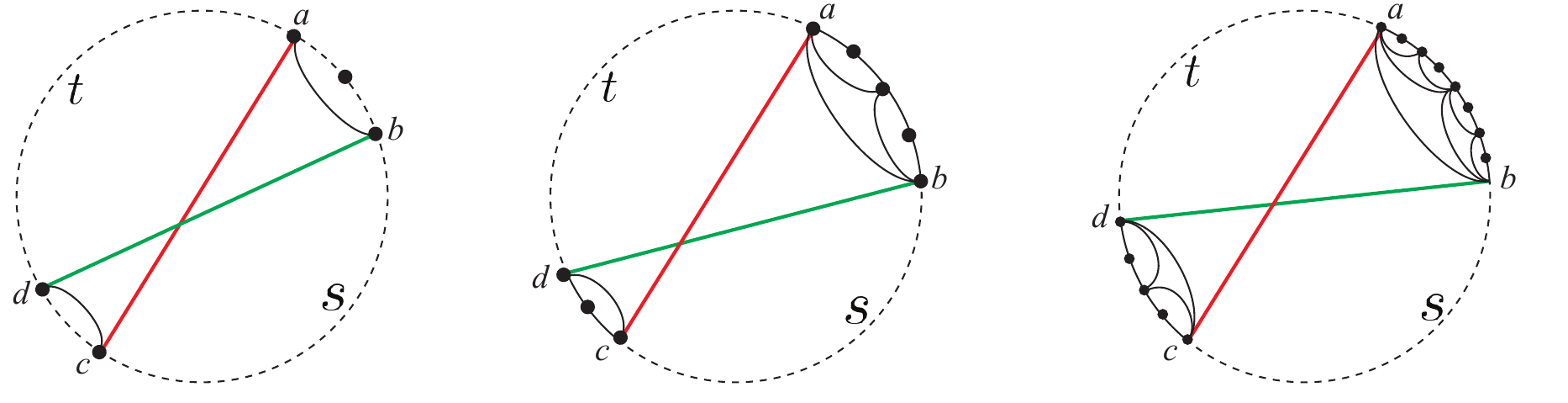}
\caption{\label{fig:replacement} Chord flip for iterations $r=1$ (left), $r=2$ (center), and $r=3$ (right)}
\end{figure}

Let $n=2^k\cdot m$, $m\ge 1$ odd and $r\le k$.
For $k\ge 3$, a flip for the iterations $r=1,2,$ and 3 of the algorithm is as on the left, center, and right of Figure \ref{fig:replacement}. 
In each figure, the chord $ac$ is replaced by the chord $bd$.
Note that on the $r$-th iteration, permanent chords of length $2^r$ are created.
On the $r$-th iteration,  if $s$ represents the length of the path included in $\mathcal{C}$ between $b$ and $c$ and $t$ represents the length of the path included in $\mathcal{C}$ between $d$ and $a$, as in Figure \ref{fig:replacement}, then a chord of length min$\{s+2^r, t+2^{r-1} \}$ is replaced by a chord of length min$\{s+2^{r-1}, t+2^r\}$. 
Such flips decrease the total chord length by up to $2^{r-1}$, creating a gap in the sequence of total chord lengths of at most $2^{r-1}$.
This gap can be remedied. We describe how this is done.
A chord of length $2^r$ in the greedy graph implies the existence of a chord of length $2^{r-1}.$
Under a chord of length $2^r$ in the greedy graph, chord changes can be executed to increase the total chord length with any value between 1 and $2^{r-1}-1$. 
In Figure \ref{fig:gaps}, under the chord $ad$, the chord $bd$ of length $2^{r-1}$ (left) is replaced by  the chord $ac$ of length $2^{r-1}+2^{r-2}$. 
This represents a $2^{r-2}$ gain in the total chord length. 
A similar flip can be made under the chord $ab$ and further, under chords incident to $a$ of lengths $2^k$, $k=r-2, r-3, \ldots, 2$,  creating gains in the total chord length of $2^{r-4}, \ldots,2, 1$, respectively.
These flips can be made independently of each other.
Since every number between 1 and $2^{r-1}-1$ can be written as a sum of numbers in the set $\{1, 2, \ldots, 2^{r-2}\}$, maximal outerplanar graphs can be generated in such a way that their total chord lengths fill the gaps in the sequence of total chord lengths for the graphs generated by the algorithm. 

For  $n=2^k\cdot m$, $m$ odd, and $r=k+1$, the $r$th iterations of the algorithm will give permanent chords of length $2^r$ but also a permanent chord of length $2^{r}+2^{r-1}$.  This is because the region of the disk that is being triangulated is bounded by an odd cycle. 
Depending on the value of $n$, odd cycles may be created on subsequent iterations and these longer chords will be layered over each other.
It is possible that on the $r$th iteration chords of length of up to  $2^{r+1}-1$ are created.   
However, a flip still decreases the total chord length by up to $2^{r-1}$, and gaps in the sequence of total chord lengths can still be filled as described above.

We conclude that between $\TCL(G_n)$ and $\TCL(F_n)$, any total chord length is possible for a graph of order $n$.\\

\begin{figure}[ht]
\centering
\includegraphics[width=0.8\textwidth]{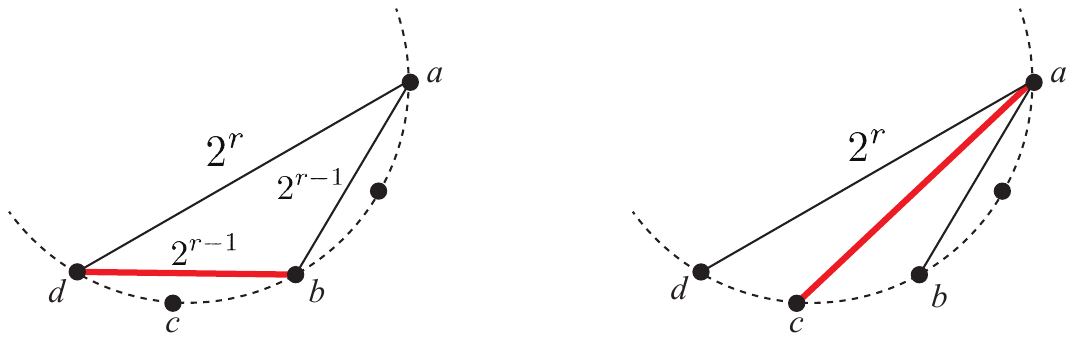}
\caption{\label{fig:gaps} Chord change under a chord of length $2^r$. The chord $bd$ (left) is replaced by the chord $ac$ (right), resulting in a total chord length gain of $2^{r-2}$.}
\end{figure}

\end{proof}
\begin{remark} The algorithm described in Theorem~\ref{thm:flip} transitions from the fan graph \(F_n\) to the greedy graph \(G_n\) through a series of re-anchoring and flipping steps. In each re-anchoring step, a new permanent chord is added, and a subset of chords incident to the previous anchor is replaced by chords incident to the new anchor. With each re-anchoring, one vertex is removed from the set of potential future anchors, and the number of chords processed in each step decreases by two. Since there are \(O(n)\) such steps, and each may require up to \(O(n)\) operations, the total time complexity of the re-anchoring phase is \(O(n^2)\).
After re-anchoring, the algorithm performs additional flips to adjust the total chord length in finer increments. 
The number of additional flips is bounded by \(O(\log n)\), since under each chord of length $2^r$ at most $r$ flips are needed. Since there are \(O(n)\) chords to choose from, and each chord is involved in at most \(O(\log n)\) flips, the time complexity for this second phase is \(O(n \log n)\). Thus, the overall time complexity of the algorithm is \(O(n^2)\), dominated by the re-anchoring phase.
\end{remark}

\noindent \textbf{Open Question:} We have found a defining characteristic for graphs with maximum TCL: they have exactly two ears. However, no analogous characteristic was found for graphs with the minimum TCL, and we leave this as an open question. 
\section{Conclusion}
\label{sec:conclusion}

In this paper, we introduced the total chord length (TCL), a new combinatorial parameter for maximal outerplanar graphs. We derived explicit formulas for both the minimum and maximum TCL values for graphs of a given order, and we gave a complete characterization of the graphs realizing the maximum value. We also showed that every integer value between the minimum and maximum TCL can be achieved by some maximal outerplanar graph of the same order. TCL is directly related to the well-known minimum-weight triangulation (MWT) problem, but it is defined in a combinatorial sense (using the vertex distance along a cycle) rather than the Euclidean length. This combinatorial perspective allows us to derive closed-form results that are difficult to obtain for the general MWT problem, which is NP-hard in the Euclidean setting. The TLC  may have applications in various domains of computational geometry and network design. For example, graphs with minimum TCL might be particularly well-suited for applications where reducing cumulative distance is desirable, such as in efficient routing or minimizing communication delays. 

\section{Acknowledgements}
The authors would like to sincerely thank the referees for going above and beyond in carefully reviewing the article. Their insightful feedback and suggestions have made the manuscript substantially stronger.
\label{sec:conclusion}

\end{document}